\documentclass[12pt,draft]{article}

\usepackage{amsfonts,amsmath,amssymb,theorem}

\usepackage[english]{babel}

\DeclareMathOperator{\dv}{div}

\DeclareMathOperator{\loc}{loc}

\DeclareMathOperator{\supp}{supp}
\DeclareMathOperator{\dist}{dist}

\begin{document}

\centerline{\bf Weighted Sobolev spaces and embedding theorems}
\vskip 0.5cm

\centerline{\bf V.~Gol'dshtein, A.~Ukhlov\footnote{The second author was partially supported by Israel Ministry of Immigrant Absorption }}
\vskip 0.5cm

\centerline{ABSTRACT.} 
\bigskip
{\small In the present paper we study embedding operators for 
weighted Sobolev spaces whose weights satisfy the well-known Muckenhoupt
$A_p$-condition. Sufficient conditions for boundedness and
compactness of the embedding operators are obtained for smooth domains and 
domains with boundary singularities. The proposed method is based on 
the concept of 'generalized' quasiconformal homeomorphisms (homeomorphisms with bounded mean distortion.) 
The choice of the homeomorphism type depends on the choice of the corresponding weighted Sobolev space. Such classes of homeomorphisms induce bounded composition operators for weighted Sobolev spaces. With the help of these homeomorphism classes the embedding problem for non-smooth domains is reduced to the corresponding classical embedding problem for smooth domains. Examples of domains with anisotropic H\"older singularities demonstrate sharpness of our machinery comparatively with known results.
} 
\vskip 0.5cm

\centerline{\bf Introduction}

\bigskip

Weighted Sobolev spaces are solution spaces of degenerate elliptic equations (see, for example [1]). Type of a weight depends on the equation type. Similar to the classical theory of Sobolev spaces, embedding theorems of weighted Sobolev spaces are suitable for the corresponding elliptic boundary problems, especially for existence and uniqueness of solutions. Embedding operators for weighted Sobolev spaces in smooth domains were studied by many authors (see, for example, [1--6]) with the help of the integral representations theory adopted to the weighted case. Weighted Sobolev spaces in non smooth domains were not studied before, except article [7] where some sufficient conditions for boundedness of the embedding operators were obtained. The main technical problem for the non-smooth case is an adequate description of an interplay between weights and boundary types (singularities). The adequate choice allows us to obtain sharp Sobolev type embeddings. 

The relation between Jacobians of quasiconformal homeomorphisms and admissible weights for Sobolev and Poincare inequalities was studied in [8].
In the present article we introduce a new approach based on the concept of 'generalized' quasiconformal homeomorphisms (or homeomorphisms with bounded mean distortion in another terminology) that induce bounded composition operators of weighted Sobolev spaces.
These homeomorphisms transform the original embedding operators on non smooth domains to the embedding operators on smooth domains with a corresponding weight change. This approach was suggested in [9] for classical Sobolev spaces on non smooth domains and can be briefly described with the help of the following diagram:
\begin{eqnarray}
W^{1}_{p}({D}',w)\overset{\varphi^{\ast}}{\longrightarrow}W^{1}_{q}(D)
\nonumber
\\
\hskip -6cm\downarrow \hskip 3cm \downarrow
\nonumber
\\
\hskip -2cm{L_s({D}',w)\overset{(\varphi^{-1})^{\ast}} {\longleftarrow}L_r(D)}
\nonumber
\end{eqnarray}
Here the operator $\varphi^{\ast}f=f\circ\varphi$ is a bounded
composition operator of Sobolev spaces induced by a homeomorphism $\varphi$ that maps smooth domain $D\subset\mathbb R^n$ onto non smooth domain $D'\subset\mathbb R^n$. Suppose that its inverse homeomorphism $\varphi^{-1}$ induces a bounded composition operator of corresponding Lebesgue spaces. 
If the Sobolev space $W^{1}_{q}(D)$ permits a bounded (compact) embedding operator into $L_r(D)$ then, using the corresponding compositions, we can construct the embedding operator of the weighted Sobolev space $W^{1}_{p}({D}',w)$ into $L_s({D}',w)$. The same scheme was used in the article [10]  for the study of the embedding operators  of $W^1_2$ into $L_2$ on non smooth bounded domains. In article [11] the same approach was applied to embedding problems for domains of Carnot groups. 

Let us shortly describe the content of the paper. In section 1 we give necessary definitions and prove density of smooth functions in weighted Sobolev spaces with weights satisfying the $A_p$-condition. Such weighted Sobolev spaces are Banach spaces.
In section 2 we introduce classes of quasiisometrical homeomorphisms and prove sufficient conditions for the compactness of the embedding operators for the weighted Sobolev spaces in domains quasiisometrically equivalent to smooth ones. In section 3 we introduce classes of homeomorphisms with bounded mean distortion and study embedding operators for weighted Sobolev spaces defined on images of smooth domains. We apply these abstract results to domains with anisotropic H\"older type singularities. The obtained estimates are sharper than the known result [7]. 

In section 5 we apply embedding theorems for weighted Sobolev spaces to degenerated elliptic boundary problems.

\bigskip

\centerline{\bf 1.~Weighted spaces}

\bigskip

In this paper we study weighted Lebesgue and Sobolev spaces defined in the domains of  $n$-dimensional
Euclidean space $\mathbb R^n$, $n\geq 2$. 

Let $D$ be an open subset of $\mathbb R^n$, $n\geq 2$, and  $w: \mathbb R^n\to [0,\infty)$ be a locally summable nonnegative function, i.e., a weight. Define a weighted Lebesgue space $L_p(D,w)$, $1\leq p<\infty$, as a Banach space of locally summable  functions $f:D\to \mathbb R$ equipped with the following norm:
$$
\|f\mid L_p(D,w)\|=
\biggr(\int\limits_D|f|^p(x)w(x)\,dx\biggr)^{1/p},\,\,\,1\leq p<\infty.
$$
Define a weighted Sobolev space
$W^m_p(D,w)$, $1\leq m<\infty$, $1\leq p<\infty$, as a normed space of locally summable, $m$ times weakly differentiable functions $f:D\to\mathbb R$ equipped with the following norm:
$$
\|f\mid W^m_p(D,w)\|=
\biggr(\int\limits_D|f|^p(x)w(x)\,dx\biggr)^{1/p}+\sum\limits_{|\alpha|=m}\biggr(\int\limits_D |D^{\alpha}f|^p(x)w(x)\,dx\biggr)^{1/p},
$$
$\alpha:=(\alpha_1,\alpha_2,...,\alpha_n)$ is a multiindex, $\alpha_i=0,1,...$, $|\alpha|=\alpha_1+\alpha_2+...+\alpha_n$ and $D^{\alpha} f$ is the weak derivative of order $\alpha$ of the function $f$:
$$
\int\limits_D f D^{\alpha}\eta~dx=(-1)^{|\alpha|}\int\limits_D (D^{\alpha}f) \eta~dx, \quad \forall \eta\in C_0^{\infty}(D).
$$ 
As usual $C_0^{\infty}(D)$ is the space of infinitely smooth functions with a compact support.  

By technical reasons we will need also a seminormed space $L^m_p(D)$ of locally summable, $m$ times weakly differentiable functions $f:D\to\mathbb R$ equipped with the following seminorm:
$$
\|f\mid L^1_p(D,w)\|=
\sum\limits_{|\alpha|=m}\biggr(\int\limits_D |D^{\alpha}f|^p(x)w(x)\,dx\biggr)^{1/p}.
$$

Without additional restrictions the space $W^m_p(D,w)$ is not necessarily a Banach space (see, for example [1]).

Let us assume additionally that
weight $w : \mathbb R^n\to [0,\infty)$ satisfies the well-known $A_p$-condition:
$$
\sup\limits_{B\subset \mathbb R^n} \biggl(\frac{1}{|B|}\int\limits_B
w\biggr) \biggl(\frac{1}{|B|}\int\limits_B
w^{1/(1-p)}\biggr)^{p-1}<+\infty,
$$
where $1<p<\infty$, and $|B|$ is Lebesgue measure of ball $B$.

{\bf Theorem~1.} Let $D\subset \mathbb R^n$ be an open set and a weight $w$ satisfies to
the $A_p$-condition. Then $W^m_p(D,w)$, $1\leq m<\infty$, $1<p<\infty$, is a Banach space. 
Smooth functions of the class $W^m_p(D,w)$ are dense in $W^m_p(D,w)$.
\vskip 0.3cm

Let us give some remarks before the proof.

Suppose that nonnegative function $\omega:\mathbb R^n\to [0,\infty)$ belongs to $C^{\infty}(\mathbb R^n)$, $\supp \omega\subset\overline{B(0,1)}$ and 
$$
\int_{\mathbb R^n}\omega~dx=1.
$$
Denote by
$$
A_{r}f(x)=\frac{1}{{r}^n}\int\limits_{\mathbb R^n} \omega\biggl(\frac{x-z}{r}\biggr)f(z)~dz
$$
a mollifier function of $f$ with a mollification kernel $\omega$.  

Let $D_{\delta} = \{x\in D :\dist(x,\partial D)>\delta\}$ for $\delta>0$. The proof of the theorem is based on  the following lemma (see, for example [12]):

\vskip 0.3cm
{\bf Lemma~1.} Let $D\subset \mathbb R^n$ be an open set and a function $f\in L_{1,\loc}(D)$ has a weak derivative $D^{\alpha}f$ on $D$. Then for every $0<r<\delta$
$$
D^{\alpha}(A_{r}f)=A_{r}(D^{\alpha}f)\quad\text{on}\quad D_{\delta}.
$$
\vskip 0.3cm

{\sc Proof.}
For the reader convenience we reproduce here a version of the proof. Note, that for every $x\in D_{\delta}$,
$$
(A_r f)(x)=\int\limits_{B(0,1)}f(x-rz)\omega(z)~dz,\quad 0<r<\delta.
$$
By  definition of the weak derivative $D_x^{\alpha}(f(x-rz))=(D_x^{\alpha}f)(x-rz))$ on $D_{\delta}$.

Suppose $(x,z)\in D_{\delta}\times B(0,1)$. Define $F(x,z)=f(x-rz)\omega(z)$ and $G(x,z)=(D^{\alpha}_xf)(x-rz)\omega(z)$. Then for each compact $K\subset D_{\delta}$ the functions $F, G$ belong to $L_1(K\times B(0,1))$. Moreover, by Fubini's theorem and the definition of the mollification kernel we get
\begin{multline}
\int\limits_K\biggl(\int\limits_{B(0,1)}|f(x-rz)\omega(z)|~dz\biggr)dx=
\int\limits_K\biggl(\frac{1}{r^n}\int\limits_{B(x,r)}|f(y)\omega\bigl(\frac{x-y}{r}\bigr)|~dy\biggr)dx\\
=\int\limits_{K^r}\biggl(\frac{1}{r^n}\int\limits_{K}|f(y)\omega\bigl(\frac{x-y}{r}\bigr)|~dx\biggr)dy\leq
\int\limits_{K^r}|f(y)|\biggl(\frac{1}{r^n}\int\limits_{K}\omega\bigl(\frac{x-y}{r}\bigr)~dx\biggr)~dy\\
\leq\int\limits_{K^r}|f(y)|\biggl(\frac{1}{r^n}\int\limits_{\mathbb R^n}\omega\bigl(\frac{x-y}{r}\bigr)~dx\biggr)~dy=
\int\limits_{K^r}|f(y)|~dy=\|f\mid L_1(K^r)\|.
\nonumber
\end{multline}
Here $K^r$ is the $r$-neighborhood of $K$, $\overline{K^r}\subset D$.
Of course, the same estimate is correct for $G$.

Using the Fubini's theorem we have 
\begin{multline}
\int\limits_{D_{\delta}}\biggl(\int\limits_{B(0,1)}D^{\alpha}_xf(x-rz)\omega(z)~dz\biggr)\eta(x)dx= \int\limits_{B(0,1)}\biggl(\int\limits_{D_{\delta}}\bigl(D^{\alpha}_xf(x-rz)\omega(z)\bigr)\eta(x)~dx\biggr)dz\\
=(-1)^{|\alpha|}\int\limits_{B(0,1)}\biggl(\int\limits_{D_{\delta}}f(x-rz)\omega(z)\bigl(D^{\alpha}\eta(x)\bigr)~dx\biggr)dz\\
=(-1)^{|\alpha|}\int\limits_{D_{\delta}}\biggl(\int\limits_{B(0,1)}f(x-rz)\omega(z)~dz\biggr)D^{\alpha}\eta(x)~dx
\nonumber
\end{multline}
for every function $\eta\in C_0^{\infty}(D_{\delta})$.

Hence for every $x\in D_{\delta}$
\begin{multline}
D^{\alpha}((A_rf)(x))=D^{\alpha}\biggl(\int\limits_{B(0,1)} f(x-rz)\omega(z)~dz\biggr)\\
=\int\limits_{B(0,1)}D^{\alpha}_x f(x-rz)\omega(z)~dz=(A_r(D^{\alpha}f))(x).
\nonumber
\end{multline}

\vskip 0.3cm

{\sc Proof of Theorem~1.}
Fix $\delta>0$. Since the weight $w$ satisfies the $A_p$-condition, Hardy-Littlewood maximal operator 
$$
Mf(x)=\sup\limits_{\delta>r>0}\frac{1}{r^n}\int\limits_{B(x,r)} f(z)~dz
$$
is bounded in $L_p(D_{\delta}, w)$ [13]. Hence
\begin{multline}
\|A_rf-f \mid L_p(D_{\delta}, w)\|=\biggl(\int\limits_{D_{\delta}} \biggl|\int\limits_{B(0,1)} f(x-rz)\omega(z)~dz-f(x)\biggr|^pw(x)~dx\biggr)^{\frac{1}{p}}\\
=\biggl(\int\limits_{D_{\delta}} \biggl|\int\limits_{B(0,1)} (f(x-rz)-f(x))\omega(z)~dz)\biggr|^pw(x)~dx\biggr)^{\frac{1}{p}}\\
\leq \|M\|\max\limits_{x\in B(0,1)}\omega(x)\biggl(\int\limits_{D_{\delta}} |f(x-rz)-f(x)|^pw(x)~dx\biggr)^{\frac{1}{p}}.
\nonumber
\end{multline}
Here $\|M\|$ is the norm of the maximal operator in the space $L_p(D_{\delta}, w)$.

From the last inequality follows, that for continuous functions $f$ 
$$
A_r f \rightarrow f \quad \text{in}\quad L_p(D_{\delta},w).
$$
Using an approximation of an arbitrary function $f\in L_p(D_{\delta},w)$  by continuous functions (see, for example, [14]) the convergence can by obtained for $f\in L_p(D_{\delta},w)$ also.

By Lemma~1 
$$
A_r f \rightarrow f \quad \text{in}\quad W^m_p(D_{\delta},w)
$$
for an arbitrary function $f\in W^{m}_{p}(D,w)$. Therefore smooth functions are dense in $W^m_p(D_{\delta},w)$.

Density of smooth functions of class $W^m_p(D,w)$ in $W^m_p(D,w)$ will be proved using the scheme proposed in [15].

Choose a sequence of open sets $D_j \Subset D_{j+1}\Subset D$, $j\geq 1$ compactly embedded one into another, such that $\cup_{j}D_j=D$. Let $\Psi$ be a partition of unity on $D$, subordinate to the covering $D_{j+1}\setminus\overline{D}_{j-1}$. Let $\psi_j$ denote the (finite) sum of those $\psi\in \Psi$ for which $\supp\psi\subset D_{j+1}\setminus\overline{D}_{j-1}$. Thus $\psi_j\in C_0^{\infty}(D_{j+1}\setminus\overline{D}_{j-1})$ and $\sum_j \psi_j\equiv 1$ in $D$. 

Fix $\varepsilon>0$ and for each $j=1,2,\dots$ choose $\varphi_j\in C_0^{\infty}(D_{j+1}\setminus\overline{D}_{j-1})$ such that 
$$
\|\varphi_j-\psi_j f \mid W^{m}_{p}(D,w)\|\leq \varepsilon 2^{-j}.
$$
Then $\varphi=\sum_j\varphi_j\in C^{\infty}(D)$, and
\begin{multline}
\|\varphi-f \mid W^{m}_{p}(D,w)\|=
\bigl\|\sum\limits_{j}\varphi_j-\sum\limits_{j}\psi_j f \mid W^{m}_{p}(D,w)\bigr\|\\
\leq \sum\limits_{j}\|\varphi_j-\psi_j f \mid W^{m}_{p}(D,w) \|<\varepsilon.
\nonumber
\end{multline}

Therefore, the weighted space
$W_{p}^{m}(D,w)$ is a Banach space and smooth functions of the class $W^{m}_{p}(D,w)$ are dense in this space.

\vskip 0.3cm

Denote by $w(A)=\int_A w(x)~dx$ the weighted measure (the measure associated with the weight $w$) of
a measurable set $A\subset\mathbb R^n$. By [13] 
the Muckenhoupt $A_p$-condition leads to the doubling condition for the weighted measure i.e ,
$\mu(B(x,2r))\leq C\alpha\mu(B(x,r))$ for all $x\in\mathbb R^n$ and
$r>0$.

{\bf Definition~1.} We call a bounded subdomain $D$ of Euclidean space $\mathbb R^n$ an
embedding domain, if for any function $f\in L^1_q(D)$, $1\leq q< n$, the Sobolev-Poincare
type inequality
$$
\inf\limits_{c\in\mathbb R}\|f-c\mid L_r(D)\|\leq M\|f\mid L^1_q(D)\|
$$
holds for any $r\leq nq/(n-q)$. Here a constant $M$ depends on $q$ and $r$ only.

Remind that for any embedding domain and for any $r < nq/(n-q)$ the corresponding embedding operator $W^1_q(D) \hookrightarrow L_r(D)$ is compact. 

Lipschitz bounded domains $D\subset\mathbb R^n$ represent examples of embedding domains. Let us remind also that 
Sobolev type embeddings for smooth domains have been studied thoroughly and a discussion about the different aspects of the embedding problem can be found, for example, in [16]. 

\bigskip

\centerline{\bf 2.~Quasiisometrical mappings and Sobolev embeddings}

\bigskip

Let $D$ and $D'$ be domains in Euclidean space $\mathbb R^n$,
$n\geq 2$. A homeomorphism $\varphi : D\to D'$
is called $Q$-quasiisometrical, if there exists a constant
$0<Q<+\infty$, such that
$$
\frac{1}{Q}\leq \underline{\varphi}^{\prime}(x)\leq \overline{\varphi}^{\prime}(x)\leq Q
$$
for all points $x \in D$.
Here
$$
\underline{\varphi}^{\prime}(x)=\liminf\limits_{z\to x}\frac{|\varphi(x)-\varphi(z)|}{|x-z|} 
\quad \text{and}\quad
\overline{\varphi}^{\prime}(x)=\limsup\limits_{z\to x}\frac{|\varphi(x)-\varphi(z)|}{|x-z|}. 
$$

It is well known that any $Q$-quasiisometrical homeomorphism is locally
bi-Lipschitz, weakly differentiable and differentiable almost
everywhere in $D$. Hence its Jacobi matrix $D\varphi =
\bigl(\frac{\partial \varphi_i}{\partial x_j}\bigr)$,
$i,j=1,\dots,n$ and its Jacobian
$J(x,\varphi)=\det\bigl(\frac{\partial \varphi_i}{\partial
x_j}\bigr)$ are well defined almost everywhere in $D$. By definition of a $Q$-quasiisometrical homeomorphism
$$
Q^{-n}\leq |J(x,\varphi)|\leq Q^n
$$
almost everywhere.

Let us remind also, that a quasiisometrical homeomorphism has Luzin
$N$-property: the image of a set of measure zero is a set of measure
zero. 

Therefore, for any $Q$-quasiisometrical homeomorphism $\varphi$ the change of variable formula in the Lebesgue integral
$$
\int\limits_{\varphi(E)}f(y)~dy=\int\limits_E f\circ\varphi(x)|J(x,\varphi)|~dx
$$
holds for any nonnegative measurable function $f$ and any measurable set $E\subset D$ [17].
 \vskip 0.3cm
 
Suppose weight $w$ satisfies the $A_p$-condition and  a homeomorphism $\varphi: \mathbb R^n \to \mathbb R^n $ is $Q$-quasi\-isometrical. Combining the change of variable formula and the two side estimate for $|J(x,\varphi)|$ we can conclude that the weight $w\circ\varphi$
also satisfies the $A_p$-condition.

\vskip 0.3cm
{\bf Theorem~2.} Let $\varphi : D\to D'$ be a $Q$-quasiisometrical homeomorphism. Then a composition operator
$\varphi^{\ast}f = f\circ\varphi$ is an isomorphism of the Sobolev spaces $W^1_p(D,w)$ and $W^1_p(D',w')$, $1< p<\infty$, where $w=w'\circ\varphi$.
\vskip 0.3cm

{\sc Proof.} Choose an arbitrary function $f \in W^1_p(D',w')$. By [14] $f$ belongs to the space $W^1_{1,\loc}(D')$  and by [17] the composition $f\circ\varphi\in W^1_{1,\loc}(D)$. Hence
\begin{multline}
\|\varphi^{\ast}f\mid W^1_p(D,w)\|=\biggl(\int\limits_D |f\circ\varphi|^p (w'\circ\varphi)(x)~dx\biggr)^{\frac{1}{p}}+
\biggl(\int\limits_D |\nabla(f\circ\varphi)|^p (w'\circ\varphi)(x)~dx\biggr)^{\frac{1}{p}}\\
\leq \biggl(\int\limits_D |f\circ\varphi|^p (w'\circ\varphi)(x)|J(x,\varphi)|\frac{1}{|J(x,\varphi)|}~dx\biggr)^{\frac{1}{p}}\\+
\biggl(\int\limits_D |\nabla f|^p(\varphi(x))(w'\circ\varphi)(x)|J(x,\varphi)|\frac{|\overline{\varphi}^{\prime}(x)|^p}{|J(x,\varphi)|}~dx\biggr)^{\frac{1}{p}}.
\nonumber
\end{multline}

Since $\varphi$ is a $Q$-quasiisometrical homeomorphism the following estimates are correct:
$$
\frac{1}{|J(x,\varphi)|}\leq Q^n \quad\text{for almost all }\quad x\in D
$$
and
$$
|\overline{\varphi}^{\prime}(x)|\leq Q \quad\text{for almost all }\quad x\in D.
$$

Hence
\begin{multline}
\|\varphi^{\ast}f\mid W^1_p(D,w)\|\leq Q^{\frac{n}{p}}
\biggl(\int\limits_D |f\circ\varphi|^p (w'\circ\varphi)(x)|J(x,\varphi)|~dx\biggr)^{\frac{1}{p}}\\+
Q^{\frac{n}{p}+1}\biggl(\int\limits_D |\nabla f|^p(\varphi(x))w'\circ\varphi(x)|J(x,\varphi)|~dx\biggr)^{\frac{1}{p}}.
\nonumber
\end{multline}

Using the change of variable formula we finally get the following inequality
\begin{multline}
\|\varphi^{\ast}f\mid W^1_p(D,w)\|\leq Q^{\frac{n}{p}}
\biggl(\int\limits_{D'} |f|^p w'(y)|~dy\biggr)^{\frac{1}{p}}+
Q^{\frac{n}{p}+1}\biggl(\int\limits_{D'} |\nabla f|^p(y)w'(y)|~dy\biggr)^{\frac{1}{p}}\\
\leq Q^{\frac{n}{p}}\bigl(Q+1\bigr)\|f\mid W^1_p(D',w')\|.
\nonumber
\end{multline}

Since the inverse homeomorphism $\varphi^{-1}$ is also a $Q$-quasiisometrical one, the  inverse inequality 
$$
\|(\varphi^{-1})^{\ast}g\mid W^1_p(D',w')\|
\leq Q^{\frac{n}{p}}\bigl(Q+1\bigr)\|g\mid W^1_p(D,w)\|, \quad g\in W^1_p(D,w).
$$
is also correct.
\vskip 0.3cm

{\bf Corollary~1.} Let $D$ and $D'$ be domains in the Euclidean space $\mathbb R^n$. Suppose there exists 
a $Q$-quasiisomerical homeomorphism $\varphi : D\to D'$.
Then the embedding operator 
$$
i : W^1_p(D,w)\hookrightarrow L_r(D,w)
$$
is bounded (compact) if and only if the embedding operator
$$
i' : W^1_p(D',w')\hookrightarrow L_r(D',w')
$$
is bounded (compact).
\vskip 0.3cm

{\sc Proof.} Suppose the embedding operator 
$$
i : W^1_p(D,w)\hookrightarrow L_r(D,w)
$$
is bounded (compact).
Since $\varphi^{-1}$ is a $Q$-quasiisometrical homeomorphism, then a composition operator
$(\varphi^{-1})^{\ast}=g\circ\varphi^{-1}:L_r(D,w) \hookrightarrow L_r(D',w')$ is bounded as one can see by a simple calculation:
\begin{multline}
\|(\varphi^{-1})^{\ast}g \mid L_r(D',w')\|\\=\biggl(\int\limits_{D'}|g\circ\varphi^{-1}|^r(y)w'(y)~dy\biggr)^{\frac{1}{r}}=\biggl(\int\limits_D |g|^r(x) w'\circ\varphi(x)|J(x,\varphi)|~dx\biggr)^{\frac{1}{r}}\\
\leq Q^{\frac{n}{r}}\biggl(\int\limits_D |g|^r(x) w(x)~dx\biggr)^{\frac{1}{r}}=Q^{\frac{n}{r}}\|g\mid L_r(D,w)\|.
\nonumber
\end{multline}
Therefore, the embedding operator $i': W^1_p(D',w')\hookrightarrow
L_r(D',w')$ is bounded (compact) as a composition of bounded operators
$\varphi^{\ast}$, $(\varphi^{-1})^{\ast}$ and bounded (compact) embedding
operator $W^1_p(D,w)\hookrightarrow L_r(D,w)$.

The proof in the inverse direction is the same.
\vskip 0.3cm

By Corollary~1 an image $D'=\varphi(D)$ of an embedding domain $D$ under a quasiisometrical homeomorphism $\varphi$ is an embedding domain also. In the paper
[10] the various examples of embedding domains of such type were discussed.

The next theorem demonstrates simple conditions for compactness
of the embedding operators of weighted Sobolev spaces.

\vskip 0.3cm
{\bf Theorem~3.} Let $D'\subset\mathbb R^n$ be a
quasiisometrical image of an embedding domain $D$, $1\leq s\leq r<nq/(n-q)$, $q\leq p$, $1<p<\infty$,
and
$$
K(w)=\max\biggl\{\|w^{-\frac{1}{p}}\mid L_{\frac{pq}{p-q}}(D')\|,\,\,\|w^{\frac{1}{s}}\mid L_{\frac{rs}{r-s}}(D')\|\biggr\}<+\infty.
$$
Then the embedding operator
$$
i: W^1_p(D',w)\hookrightarrow L_s(D',w)
$$
is a compact operator.

For $r=nq/(n-q)$ the embedding operator $i$ is bounded only.
\vskip 0.3cm

{\sc Proof.} By conditions of the theorem there
exists a $Q$-quasiisometrical homeomorphism $\varphi : D\to D'$ of the embedding domain $D$ onto the domain $D'$. For any
function $u\in W^1_p(D',w)$ the composition $u\circ \varphi$ is
weakly differentiable in the domain $D$, and the following
estimate is correct
\begin{multline}
\|u\circ\varphi\mid W^1_q(D)\|=\biggl(\int\limits_D
|u\circ\varphi|^q~dx\biggr)^{\frac{1}{q}}+\biggl(\int\limits_D
|\nabla (u\circ\varphi)|^q~dx\biggr)^{\frac{1}{q}}\\
\leq
\biggl(\int\limits_D |u\circ\varphi|^q
(|J(x,\varphi)|w(\varphi(x))^{\frac{q}{p}}\frac{1}{(|J(x,\varphi)|w(\varphi(x))^{\frac{q}{p}}}
~dx\biggr)^{\frac{1}{q}}\\
+\biggl(\int\limits_D |\nabla u|^q
|\overline{\varphi}^{\prime}(x)|^q(|J(x,\varphi)|w(\varphi(x))^{\frac{q}{p}}\frac{1}{(|J(x,\varphi)|w(\varphi(x))^{\frac{q}{p}}}~dx\biggr)^{\frac{1}{q}}.
\nonumber
\end{multline}

By the H\"older inequality 
\begin{multline}
\|u\circ\varphi\mid W^1_q(D)\|\leq\\
\leq
\biggl(\int\limits_D\biggl(\frac{1}{|J(x,\varphi)|w(\varphi(x))}\biggl)^{\frac{q}{p-q}}~dx\biggr)^{\frac{p-q}{pq}}
\biggl(\int\limits_D
|u|^p(\varphi(x))w(\varphi(x))|J(x,\varphi)|~dx\biggr)^{\frac{1}{p}}\\
+\biggl(\int\limits_D\biggl(\frac{|\overline{\varphi}^{\prime}(x)|^p}{|J(x,\varphi)|w(\varphi(x))}\biggl)^{\frac{q}{p-q}}~dx\biggr)^{\frac{p-q}{pq}}
\biggl(\int\limits_D |\nabla
u|^p(\varphi(x))w(\varphi(x))|J(x,\varphi)|~dx\biggr)^{\frac{1}{p}}.
\nonumber
\end{multline}

Since $\varphi$ is a $Q$-quasiisometrical homeomorphism, then
by the change of variable formula for the Lebesgue integral
we obtain
\begin{multline}
\|u\circ\varphi\mid W^1_q(D)\|\leq Q^{\frac{n}{q}}
\biggl(\int\limits_{D'}w(y)^{\frac{q}{q-p}}~dy\biggr)^{\frac{p-q}{pq}}\biggl(\int\limits_{D'}|u|^p
w(y) ~dy\biggr)^{\frac{1}{p}}\\
+Q^{\frac{p-q+n}{q}}
\biggl(\int\limits_{D'}w(y)^{\frac{q}{q-p}}~dy\biggr)^{\frac{p-q}{pq}}\biggl(\int\limits_{D'}|\nabla
u|^p w(y)~dy\biggr)^{\frac{1}{p}}\\
\leq Q^{\frac{n}{q}}K(w)\|u\mid
L_p(D',w)\|+Q^{\frac{p-q+n}{q}}K(w)\|\nabla u\mid L_p(D'w)\|.
\nonumber
\end{multline}

By the previous inequality the composition operator 
$$
\varphi^{\ast}: W^1_p(D',w)\to W^1_q(D),\,\,\, 1\leq q\leq
p<+\infty,
$$
is bounded.

Let us prove boundedness of the composition operator $(\varphi^{-1})^{\ast}: L_r(D)\to L_s(D',w)$.
By the theorem's conditions the quantity
$\|w^{\frac{1}{s}}\mid L_{\frac{rs}{r-s}}(D')\|$ is finite. Hence by [18] the 
composition operator 
$$
(\varphi^{-1})^{\ast}: L_r(D)\to L_s(D',w), \,\,\, 1\leq s\leq
r<+\infty,
$$
is bounded.

Because $D$ is an embedding domain the embedding operator $i:W^1_q(D)\hookrightarrow L_r(D)$ is compact for any $r<nq/(n-q)$ and bounded for $r=nq/(n-q)$.
Therefore the embedding operator $i': W^1_p(D',w)\hookrightarrow
L_s(D',w)$ is compact (bounded) as a composition of bounded operators
$\varphi^{\ast}$, $(\varphi^{-1})^{\ast}$ and the compact (bounded) embedding
operator $i$ for any 
$r<nq/(n-q)$ ($r=nq/(n-q)$).
\vskip 0.3cm

In a similar way we can prove

\vskip 0.3cm
{\bf Theorem~4.} Let $D'\subset\mathbb R^n$ be a
quasiisometrical image of an embedding domain $D$, $1\leq s\leq r<nq/(n-q)$, $q\leq p$, $1<p<\infty$,
and
$$
K(w)=\|w^{\frac{1}{s}}\mid L_{\frac{rs}{r-s}}(D')\|<+\infty.
$$
Then the embedding operator
$$
i: W^1_p(D')\hookrightarrow L_s(D',w)
$$
is compact.

For $r=nq/(n-q)$ the embedding operator $i$ is bounded only.
\vskip 0.3cm

The next lemma allows us to construct various examples of the embedding domains.

\vskip 0.3cm
{\bf Lemma~2.} Let $D_1$ and $D_2$ be such domains that embedding operators 
$$
i: W^1_p(D_1, w)\hookrightarrow L_s(D_1,v),
$$
$$
i: W^1_p(D_2, w)\hookrightarrow L_s(D_2,v)
$$
are compact.
Then the embedding operator 
$$
i: W^1_p(D_1\cup D_2, w)\hookrightarrow L_s(D_1\cup D_2,v)
$$
is also compact.

\vskip 0.3cm

{\sc Proof.} We prove this lemma by the scheme suggested in [10].
Choose a sequence of functions $\{f_n\}\subset W^1_p(D_1\cup D_2, w)$ such that $\|f_n\mid W^1_p(D_1\cup D_2, w)\|\leq 1$ for all $n$.
Let $g_n=f_n\vert_{D_1}$ and $h_n=f_n\vert_{D_2}$. Then $g_n\in W^1_p(D_1)$, $h_n\in W^1_p(D_2)$, $\|f_n\mid W^1_p(D_1)\|\leq 1$, $\|h_n\mid W^1_p(D_2)\|\leq 1$.

Because the embedding operator $i: W^1_p(D_1,w)\hookrightarrow L_s(D_1,v)$ is compact, we can choose a subsequence $\{g_{n_k}\}$ of the sequence $\{g_n\}$ which converges in $L_s(D_1,v)$ to a function $g_0\in L_s(D_1,v)$. Because the second embedding operator $i: W^1_p(D_2,w)\hookrightarrow L_s(D_2,v)$ is also compact we can choose a subsequence $\{h_{n_{k_m}}\}$ of the sequence $\{h_{n_{k}}\}$ which converges in $L_s(D_2,v)$ to a function $h_0\in L_s(D_2,v)$. It is evident that $g_0=h_0$ $v$-almost everywhere in $D_1\cap D_2$ and the function $f_0=g_0$ on $D_1$ and $f_0=h_0$ on $D_2$ belongs to $L_s(D_1\cup D_2,v)$.

Hence
$$
\|f_{n_{k_m}}-f_0\mid L_s(D_1\cup D_2,v)\|\\
\leq\|g_{n_{k_m}}-g_0\mid L_s(D_1,v)\|+\|h_{n_{k_m}}-h_0\mid L_s(D_2,v)\|.
$$

Therefore
$\|f_{n_{k_m}}-f_0\mid L_s(D_1\cup D_2,v)\|\longrightarrow 0$ for $m\longrightarrow\infty$.

\vskip 0.3cm

\bigskip

\centerline{\bf 3.~Embedding operators for general domains}

\bigskip

Let $D$ and $D'$ be domains in Euclidean space $\mathbb R^n$,
$n\geq 2$. Remember that a homeomorphism $\varphi:D \to D'$ belongs to Sobolev class $W^1_{1,\loc}(D)$ if its coordinate functions belong to $W^1_{1,\loc}(D)$. Denote by $D\varphi$ the weak differential of $\varphi$. The norm $|D\varphi(x)|$ is the standard norm of the linear operator defined by  $D\varphi(x)$. 

Call a homeomorphism $\varphi: D\to D'$ $w$-weighted $(p,q)$-quasiconformal, if $\varphi$ belongs
to the Sobolev space $W^1_{1,\loc}(D)$, $|D\varphi|=0$ almost everywhere on the set $Z=\{x: |J(x,\varphi)|w(\varphi(x))=0\}$ and the following inequality
$$
K_{p,q}(D,w)=\biggl[\int\limits_{D}\biggl(\frac{|D\varphi(x)|^p}{|J(x,\varphi)|w(\varphi(x))}\biggr)^{\frac{q}{p-q}}~dx\biggl]^{\frac{p-q}{pq}}<\infty
$$
is correct.

For $w\equiv 1$ we call a $1$-weighted $(p,q)$-quasiconformal homeomorphism a $(p,q)$-quasi\-con\-for\-mal one.

The following result was proved in [19] for a more general class of mappings. For readers 
convenience we reproduce here a simple version of the proof adopted to homeomorphisms.

\vskip 0.3cm
{\bf Proposition~1.} Let $D$ and $D'$ be domains in Euclidean space $\mathbb R^n$,
$n\geq 2$ and $\varphi: D\to D'$ be a $w$-weighted $(p,q)$-quasiconformal homeomorphism. Then 
a composition operator
$$
\varphi^{\ast}: L^1_p(D',w)\to L^1_q(D), \,\,\, 1\leq q\leq
p<+\infty,
$$
is bounded.
\vskip 0.3cm

{\sc Proof.} Let $f \in L^1_p(D',w)$ be a smooth function. Then $f\circ\varphi \in L^1_{1,\loc}(D)$ and the following inequalities
\begin{multline}
\|\varphi^{\ast}f\mid L^1_q(D)\|= \biggl(\int\limits_D |\nabla(f\circ\varphi)|^q~dx\biggr)^{\frac{1}{q}}
\leq  \biggl(\int\limits_D |D\varphi|^q |\nabla f|^q ~dx\biggr)^{\frac{1}{q}}\\
\leq \biggl(\int\limits_D |D\varphi|^q\frac{1}{|J(x,\varphi)|^{\frac{q}{p}}w(\varphi(x))^{\frac{q}{p}}} |\nabla f|^q |J(x,\varphi)|^{\frac{q}{p}}w(\varphi(x))^{\frac{q}{p}}~dx\biggr)^{\frac{1}{q}}
\nonumber
\end{multline}
are correct.

Using H\"older inequality and the change of variable formula for Lebesgue integral, we obtain
\begin{multline}
\|\varphi^{\ast}f\mid L^1_q(D)\|\leq \biggl[\int\limits_{D}\biggl(\frac{|D\varphi(x)|^p}{|J(x,\varphi)|w(\varphi(x))}\biggr)^{\frac{q}{p-q}}~dx\biggl]^{\frac{p-q}{pq}}\biggl(\int\limits_D |\nabla f|^p |J(x,\varphi)|w(\varphi(x))~dx\biggr)^{\frac{1}{p}}\\
=K_{p,q}(D,w)\|f\mid L^1_p(D',w)\|.
\nonumber
\end{multline}
The fulfillment of the last inequality for an arbitrary function $f\in L^1_p(D',w)$ can be proved by an approximation of $f$ by smooth functions [19].

\vskip 0.3cm

The next theorem gives sufficient condition for boundedness (compactness) embedding operators in non-smooth domains.

\vskip 0.3cm
{\bf Theorem~5.} Let a domain $D\subset\mathbb R^n$ be an embedding
domain  and there exists a $w$-weighted
$(p,q)$-quasiconformal homeomorphism $\varphi: D\to D'$ of domain $D$ onto bounded domain $D'$. 

If for some $p\leq s\leq r<\infty$ the following inequality is correct
$$
\int\limits_D\biggl(|J(x,\varphi)|w(\varphi(x))\biggr)^{\frac{r}{r-s}}~dx<+\infty,
$$
then an embedding operator
$$
i: W^1_p(D',w)\hookrightarrow L_s(D',w)
$$
is bounded, if $s\leq r\leq nq/(n-q)$, and is compact if $s\leq
r<nq/(n-q)$.
\vskip 0.3cm

{\sc Proof.} Because 
$$
\int\limits_D\biggl(|J(x,\varphi)|w(\varphi(x))\biggr)^{\frac{r}{r-s}}<+\infty,
$$
the composition operator $(\varphi^{-1})^{\ast}:L_r(D)\to L_s(D',w)$  is bounded, i.e. the following inequality
$$
\|(\varphi^{-1})^{\ast}v\mid L_s(D',w)\|\leq A_{r,s}(D,w) \|v\mid
L_r(D)\|
$$
is correct. Here $A_{r,s}(D,w)$ is a positive constant.

Because domain $D$ is an embedding domain and the composition operators 
$$
(\varphi^{-1})^{\ast}:L_r(D)\to L_s(D',w),\,\, \varphi^{\ast}: L^1_p(D',w)\to L^1_q(D)
$$
are bounded, the following inequalities
\begin{multline} 
\inf\limits_{c\in\mathbb R}\|u-c\mid L_s(D',w)\|\leq
A_{r,s}(D,w)\inf\limits_{c\in\mathbb R}\|v-c\mid L_r(D)\|\\ \leq
A_{r,s}(D,w)M\|v\mid L^1_q(D)\|\leq
A_{r,s}(D,w)K_{p,q}(D,w)M\|u\mid L^1_p(D',w)\|
\nonumber
\end{multline}
hold. Here $M$ and $K_{p,q}(D,w)$ are positive constants.

The H\"older inequality implies the following estimate
\begin{multline} |c|=w(D')^{-\frac{1}{p}}\|c\mid L_p(D',w)\|\leq
w(D')^{-\frac{1}{p}}\bigl(\|u\mid L_p(D',w)\|+\|u-c\mid
L_p(D',w)\|\bigr)\\
\leq w(D')^{-\frac{1}{p}}\|u\mid
L_p(D',w)\|+w(D')^{-\frac{1}{s}}\|u-c\mid L_s(D',w)\|.
\nonumber
\end{multline}

Because $q\leq r$ we have
\begin{multline} \|v\mid L_q(D)\|\leq \|c\mid L_q(D)\|+\|v-c\mid
L_q(D)\|\leq
|c||D|^{\frac{1}{q}}+|D|^{\frac{r-q}{r}}\|v-c\mid L_r(D)\|\\
\leq\biggl(w(D')^{-\frac{1}{p}}\|u\mid
L_p(D',w)\|+w(D')^{-\frac{1}{s}}\|u-c\mid
L_s(D',w)\|\biggr)|D|^{\frac{1}{q}}+|D|^{\frac{r-q}{r}}\|v-c\mid
L_r(D)\|.
\nonumber
\end{multline}

From previous inequalities we obtain finally
\begin{multline}
\|v\mid L_q(D)\|\leq
|D|^{\frac{1}{q}}w(D')^{-\frac{1}{p}}\|u\mid
L_p(D')\|\\
+A_{r,s}(D,w)K_{p,q}(D,w)M|D|^{\frac{1}{q}}w(D')^{-\frac{1}{p}}\|u\mid
L^1_p(D',w)\|\\
+K_{p,q}(D,w)M|D|^{\frac{r-q}{r}}\|u\mid L^1_p(D',w)\|.
\nonumber
\end{multline}
Therefore the composition operator  
$$
\varphi^{\ast}: W^1_p(D',w)\to W^1_q(D)
$$
is also bounded.

Finally we can conclude that the embedding operator $i: W^1_p(D',w)\hookrightarrow
L_s(D',w)$ is bounded as the composition of bounded operators
$\varphi^{\ast}$, $(\varphi^{-1})^{\ast}$ and the embedding operator
$W^1_q(D)\hookrightarrow L_r(D)$ in the case $r\leq nq/(n-q)$.
The embedding operator $i: W^1_p(D',w)\hookrightarrow
L_s(D',w)$ is compact as the composition of bounded operators
$\varphi^{\ast}$, $(\varphi^{-1})^{\ast}$ and compact the embedding
operator $W^1_q(D)\hookrightarrow L_r(D)$, in the case
$r<nq/(n-q)$.
\vskip 0.3cm

Let us apply these general results to domains with
anisotropic H\"older singularities described in [9]. 

Let $g_i(\tau)=\tau^{\gamma_i}$, $\gamma_i\geq 1$, $0\leq\tau\leq 1$. For the function $G=\prod_{i=1}^{n-1}g_i$ denote by
$$
\gamma=\frac{\log G(\tau)}{\log \tau}+1.
$$
It is evident that $\gamma\geq n$. Let us consider the domain
$$
H_g=\{ x\in\mathbb R^n : 0<x_n<1, 0<x_i<g_i(x_n),
\,i=1,2,\dots,n-1\}.
$$
In the case $g_1=g_2=\dots=g_{n-1}$ we will say that domain $H_g$ is a domain with $\sigma$-H\"older singularity, $\sigma=(\gamma-1)/(n-1)$.
For $g_1(\tau)=g_2(\tau)=\dots=g_{n-1}(\tau)=\tau$ we will use notation $H_1$ instead
of $H_g$. The domain $H_1$ is quasiisometrically homeomorphic to the
standard unit ball. Hence by Theorem 3 domain $H_1$ is an embedding domain. 

If the weight is polynomial i.e. $w(x):=|x|^{\alpha}$ then the $A_p$-condition is correct only if $-n<\alpha<n(p-1)$.

\vskip 0.3cm

{\bf Theorem~6.} Let $-n<\alpha<n(p-1)$ and  $1<p<\alpha+\gamma$, then the embedding
operator
$$W^1_p(H_g,|x|^{\alpha})\hookrightarrow
L_s(H_g,|x|^{\alpha})
$$
is compact for any $1\leq
s<\frac{(\alpha+\gamma)p}{\alpha+\gamma-p}$.
\vskip 0.3cm

{\sc Proof.} For any $0<a<1$ we define a homeomorphism
$\varphi_a : H_1\to H_g$, $a>0$, by the expression
$$
\varphi_a(x)=(\frac{x_1}{x_n}g^a_1(x_n),\dots,\frac{x_{n-1}}{x_n}g^a_{n-1}(x_n),x_n^a).
$$
During the proof we will choice a number $a$ and a corresponding homeomorphism in such a way that the conditions of Theorem 5 will be fulfilled.
By simple calculation 
$$
\frac{\partial(\varphi_a)_i}{\partial
x_i}=\frac{g^a_i(x_n)}{x_n},\quad
\frac{\partial(\varphi_a)_i}{\partial
x_n}=\frac{-x_ig^a_i(x_n)}{x_n^{2}}+\frac{ax_ig^{a-1}_i(x_n)}{x_n}g'_i(x_n)
\quad\text{and}\quad\frac{\partial(\varphi_a)_n}{\partial
x_n}=ax_n^{a-1}
$$
for any $i=1,...,n-1$.

Hence $J(x,\varphi_a(x))=ax_n^{a-n}G^a(x_n)$. By definition functions $g_i$,
$i=1,2,\dots,n-1$ are Lipschitz functions. Therefore there exists
a constant $M<+\infty$ such that
$$
g_i(x_n)\leq Mx_n \quad\text{and}\quad g'_i(x_n)\leq M
$$
for any $x_n\in [0,1]$ and $i=1,2,\dots,n-1$. Using estimates for derivatives and the inequalities $x_i\leq x_n$ that are correct for all $x\in H_1$ we obtain the following estimate $|D\varphi_a(x)|\leq c_1
x_n^{a-1}$. By the same way we obtain also the two-sided estimate:
$$
c_2x_n^{\alpha a}\leq |\varphi_a(x)|^{\alpha}\leq c_3 x_n^{\alpha
a}.
$$
Now we can check for which $q$ the homeomorphism $\varphi_a : H_1\to H_g$ is a $w$-weighted $(p,q)$-quasiconformal homeomorphism:

\begin{multline}
I_a=K_{p,q}(H_1,w)^{\frac{pq}{p-q}}=\int\limits_{H_1}\biggl(\frac{|D\varphi_a(x)|^p}{|J(x,\varphi_a)|w(\varphi_a(x))}\biggr)^{\frac{q}{p-q}}~dx\\
\leq
C\int\limits_0^1\int\limits_0^{x_n}\dots\int\limits_0^{x_n}\biggl(\frac{x_n^{p(a-1)}}{x_n^{a-n}G^a(x_n)x_n^{\alpha
a }}\biggr)^{\frac{q}{p-q}}~dx_1\dots dx_{n-1}dx_n\\
=C\int\limits_0^1
x_n^{\bigl(p(a-1)-a(\alpha+1)+n\bigr)\frac{q}{p-q}+n-1}G^{-a\frac{q}{p-q}}(x_n)~dx_n.
\nonumber
\end{multline}
Hence the quantity $I_a$ is finite if
$$
\bigl(p(a-1)-a(\alpha+1)+n\bigr)\frac{q}{p-q}+n-a(\gamma-1)\frac{q}{p-q}>0
$$
or
\begin{equation}
q<np/\bigr(a(\alpha+\gamma)+p-ap\bigl). 
\nonumber
\end{equation}
Hence, the homeomorphism $\varphi_a$ is a $w$-weighted $(p,q)$-quasiconformal homeomorphism.

Let us check the conditions of Theorem 5. First we have to estimate the degree of integrability for Jacobian $J_a$ of the homeomorphism
$\varphi_a$:
\begin{multline}
J_a=\int\limits_{H_1}\biggl(|J(x,\varphi)|w(\varphi(x))\biggr)^{\frac{r}{r-s}} dx\\
\leq
C\int\limits_0^1\int\limits_0^{x_n}\dots\int\limits_0^{x_n}x_n^{\bigl(a(\alpha+1)-n\bigr)\frac{r}{r-s}}G^{a\frac{r}{r-s}}(x_n)~dx_1\dots
dx_{n-1}dx_n\\
\leq C\int\limits_0^1x_n^{\bigl(a(\alpha+1)-n\bigr)\frac{r}{r-s} +n-1
+a\frac{r}{r-s}(\gamma-1)}~dx_n.
\nonumber
\end{multline}
The integral $J_a$ converges, if
$$
\bigl(a(\alpha+1)-n\bigr)\frac{r}{r-s} +n+a\frac{r}{r-s}(\gamma-1)>0,
$$
or
$$
s<\frac{a(\alpha+\gamma)}{n}r.
$$

Hence, the conditions of Theorem~5 are fulfilled if
$$
s<\frac{a(\alpha+\gamma)}{n}r,\,\,\, r<\frac{nq}{n-q}\,\,\,\text{and}\,\,\,q<\frac{np}{a(\alpha+\gamma)+p-ap}. 
$$
Therefore
$$
s<\frac{a(\alpha+\gamma)}{n}\frac{nq}{n-q}<\frac{np}{a(\alpha+\gamma-p)}\frac{a(\alpha+\gamma)}{n}=
\frac{p(\alpha+\gamma)}{\alpha+\gamma-p}.
$$
Theorem~6 is proved.

\vskip 0.3cm

{\bf Remark~1.} The conclusion of Theorem~6 are fulfilled for functions 
$g_i:[0,1]\to \mathbb R$, $i=1,2,\dots,n-1$ such that 
$$
C_1\tau^{\gamma_i}\leq g_i(\tau)\leq C_2\tau^{\gamma_i}
$$
for some constants $C_1$ and $C_2$.
\vskip 0.3cm

From theorem~6 immediately follows.
\vskip 0.3cm

{\bf Corollary~2.} Let $D\subset\mathbb R^n$ be a domain with $\sigma$-H\"older singularity. Then the embedding operator
$$W^1_p(D,|x|^{\alpha})\hookrightarrow
L_{s}(D,|x|^{\alpha})
$$
is compact for 
$$
s\leq
\frac{(\sigma(n-1)+1+\alpha)p}{\sigma(n-1)+\alpha-(p-1)}.
$$
(Here $s\geq 0$, since $p<\alpha+\gamma$, $\sigma=(\gamma-1)/(n-1)$.)
\vskip 0.3cm

Let us compare this result with the known ones.
From the main result of [7] it follows that for an arbitrary domain $D$ with $\sigma$-H\"older singularity the embedding operator
$$W^1_p(D,|x|^{\alpha})\hookrightarrow
L_{\tilde{s}}(D,|x|^{\alpha})
$$
is bounded, while 
$$
\tilde{s}\leq \frac{(n+\alpha)p}{\sigma (\alpha
+ n-1)-(p-1)}.
$$

Then
$s>\tilde{s}$ while $\sigma>1$, and $s=\tilde{s}$ while
$\sigma=1$. Hence our estimate is sharper. 
\vskip 0.3cm

The next results deal with embeddings of classical Sobolev spaces into weighted \newline Lebesgue spaces.

{\bf Theorem~7.} Let $D\subset\mathbb R^n$ be an embedding 
domain and  there exists a
$(p,q)$-quasi\-con\-for\-mal homeomorphism $\varphi$ of $ D$  onto bounded domain $D'$. 

If 
$$
\int\limits_D\biggl(|J(x,\varphi)|w(\varphi(x))\biggr)^{\frac{r}{r-s}}<+\infty,
$$
for a pair of numbers $p\leq s<r<\infty$
then an embedding operator
$$
i: W^1_p(D')\hookrightarrow L_s(D',w)
$$
is compact for $ r < nq/(n-q)$ and is bounded for $r=nq/(n-q)$.
\vskip 0.3cm

{\sc Proof.} Because $D$ is an embedding domain $(p,q)$-quasiconformal homeomorphism $\varphi$ induces the bounded composition operator
$$
\varphi^{\ast} : W^1_p(D')\to W^1_q(D)
$$
(see [18]). 

Because $
\int\limits_D\biggl(|J(x,\varphi)|w(\varphi(x))\biggr)^{\frac{r}{r-s}}<+\infty
$  the composition operator for Lebesgue
spaces
$$
\|(\varphi^{-1})^{\ast}v\mid L_s(D',w)\|\leq A_{r,s}(D,w) \|v\mid
L_r(D)\|
$$
is bounded also.

Finally we can conclude that 

1) If $r\leq nq/(n-q)$ then the embedding operator $i: W^1_p(D')\hookrightarrow
L_s(D',w)$ is bounded as a composition of bounded operators
$\varphi^{\ast}$, $(\varphi^{-1})^{\ast}$ and the bounded embedding operator
$W^1_q(D)\hookrightarrow L_r(D)$ ;

2) If $r<nq/(n-q)$ then the embedding operator  $i: W^1_p(D')\hookrightarrow
L_s(D',w)$ is compact as a composition of bounded operators
$\varphi^{\ast}$, $(\varphi^{-1})^{\ast}$ and the compact embedding
operator $W^1_q(D)\hookrightarrow L_r(D)$.
\vskip 0.3cm

Apply the previous result to anisotropic H\"older domains.
\vskip 0.3cm

{\bf Theorem~8.} Let $1<p<\gamma$ and $1\leq
s<\frac{(\alpha+\gamma)p}{\gamma-p}$, then the embedding operator
$$W^1_p(H_g)\hookrightarrow
L_s(H_g,|x|^{\alpha})
$$
is compact.
\vskip 0.3cm

{\sc Proof.} Similar to the proof of Theorem 6 for any $0<a<1$ we define a homeomorphism
$\varphi_a : H_1\to H_g$, $a>0$ by the expression
$$
\varphi_a(x)=(\frac{x_1}{x_n}g^a_1(x_n),\dots,\frac{x_{n-1}}{x_n}g^a_{n-1}(x_n),x_n^a).
$$
During the proof we will choice a number $a$ and the corresponding homeomorphism in such a way that the conditions of Theorem 7 will be fulfilled.

By simple calculation we have
$$
\frac{\partial(\varphi_a)_i}{\partial
x_i}=\frac{g^a_i(x_n)}{x_n},\quad
\frac{\partial(\varphi_a)_i}{\partial
x_n}=\frac{-x_ig^a_i(x_n)}{x_n^{2}}+\frac{ax_ig^{a-1}_i(x_n)}{x_n}g'_i(x_n)
\quad\text{and}\quad\frac{\partial(\varphi_a)_n}{\partial
x_n}=ax_n^{a-1}
$$
for any $i=1,...,n-1$.

Hence $J(x,\varphi_a(x)=ax_n^{a-n}G^a(x_n)$. By definition functions $g_i$,
$i=1,2,\dots,n-1$ are Lipschitz. So there exists
a constant $M<+\infty$ such that
$$
g_i(x_n)\leq Mx_n \quad\text{and}\quad g'_i(x_n)\leq M
$$
for any $x_n\in [0,1]$ and $i=1,2,\dots,n-1$. Using estimates for derivatives and the inequalities $x_i\leq x_n$ that are correct for all $x\in H_1$ we obtain the following estimate $|D\varphi_a(x)|\leq c_1
x_n^{a-1}$. In the same way we obtain also the two-sided estimate
$$
c_2x_n^{\alpha a}\leq |\varphi_a(x)|^{\alpha}\leq c_3 x_n^{\alpha
a}.
$$
Now we can check for which $q$ the homeomorphism $\varphi_a : H_1\to H_g$ is a $(p,q)$-quasicon\-for\-mal homeomorphism. Let us start from the following estimate

\begin{multline}
I_a=K_{p,q}(H_1)^{\frac{pq}{p-q}}=\int\limits_{H_1}\biggl(\frac{|D\varphi_a(x)|^p}{|J(x,\varphi_a)|}\biggr)^{\frac{q}{p-q}}~dx\\
\leq
C\int\limits_0^1\int\limits_0^{x_n}\dots\int\limits_0^{x_n}\biggl(\frac{x_n^{p(a-1)}}{x_n^{a-n}G^a(x_n)}\biggr)^{\frac{q}{p-q}}~dx_1\dots dx_{n-1}dx_n\\
=C\int\limits_0^1
x_n^{\bigl(p(a-1)-a+n\bigr)\frac{q}{p-q}+n-1}G^{-a\frac{q}{p-q}}(x_n)~dx_n.
\nonumber
\end{multline}
Hence, the quantity $I_a$ is finite if
$$
\bigl(p(a-1)-a+n\bigr)\frac{q}{p-q}+n-a(\gamma-1)\frac{q}{p-q}>0
$$
or
\begin{equation}
q<np/\bigr(a\gamma+p-ap\bigl). 
\nonumber
\end{equation}
Hence, the homeomorphism $\varphi_a$ is a $(p,q)$-quasiconformal homeomorphism.

Let us check the conditions of Theorem 5. First we have to estimate the degree of integrability for Jacobian $J_a$ of the homeomorphism
$\varphi_a$:
\begin{multline}
J_a=\int\limits_{H_1}\biggl(|J(x,\varphi)|w(\varphi(x))\biggr)^{\frac{r}{r-s}} dx\\
\leq
C\int\limits_0^1\int\limits_0^{x_n}\dots\int\limits_0^{x_n}x_n^{\bigl(a(\alpha+1)-n\bigr)\frac{r}{r-s}}G^{a\frac{r}{r-s}}(x_n)~dx_1\dots
dx_{n-1}dx_n\\
\leq C\int\limits_0^1x_n^{\bigl(a(\alpha+1)-n\bigr)\frac{r}{r-s} +n-1
+a\frac{r}{r-s}(\gamma-1)}~dx_n.
\nonumber
\end{multline}
The integral $J_a$ converges, if
$$
\bigl(a(\alpha+1)-n\bigr)\frac{r}{r-s} +n+a\frac{r}{r-s}(\gamma-1)>0,
$$
or
$$
s<\frac{a(\alpha+\gamma)}{n}r.
$$

Hence, the conditions of Theorem~5 are fulfilled if
$$
s<\frac{a(\alpha+\gamma)}{n}r,\,\,\, r<\frac{nq}{n-q}\,\,\,\text{and}\,\,\,q<\frac{np}{a\gamma+p-ap}. 
$$
Therefore
$$
s<\frac{a(\alpha+\gamma)}{n}\frac{nq}{n-q}<\frac{np}{a(\gamma-p)}\frac{a(\alpha+\gamma)}{n}=
\frac{p(\alpha+\gamma)}{\gamma-p}.
$$

\bigskip

\centerline{\bf 4.~Sobolev embeddings for spaces with high derivatives}

\bigskip

This section is devoted to embedding theorems for weighted Sobolev spaces with high derivatives.

If for some $m$ and $p$ an embedding theorem for classical Sobolev spaces $W^m_p$ is correct then, using a standard procedure, it is possible to obtain the corresponding embedding theorem for  $m_1\geq m$ and $p_1\geq p$ also (see, for example [9]). Here we adopt the scheme of [9] to the case of weighted Sobolev spaces. The next lemma is the main technical result of this section.
\vskip 0.3cm

{\bf Lemma~3.} Let $D$ be a domain in $\mathbb R^n$. Suppose that for some $p_0\geq 1$ and $q_0\geq 1$ the embedding operator 
$$
W^1_{p_0}(D,w) \hookrightarrow L_{q_0}(D,w)
$$
is bounded.
Let $p\geq p_0$ and $\frac{1}{p}>\frac{1}{p_0}-\frac{1}{q_0}$. If $q$ is such that
$$
\frac{1}{p}-\frac{1}{q}=\frac{1}{p_0}-\frac{1}{q_0}
$$
then the embedding operator 
$$
W^1_{p}(D,w) \hookrightarrow L_{q}(D,w)
$$
is bounded also.
\vskip 0.3cm

{\sc Proof.} Let a function $u$ belongs to the space $W^1_p(D,w)\cap C_0^{\infty}(D)$. Using boundedness of the embedding operator
$$
W^1_{p_0}(D,w) \hookrightarrow L_{q_0}(D,w)
$$
we obtain the following estimate
\begin{multline}
\biggl(\int\limits_D |u|^q w(x)~dx\biggr)^{\frac{1}{q_0}}=\biggl(\int\limits_D \bigl(|u|^{\frac{q}{q_0}}\bigr)^{q_0} w(x)~dx\biggr)^{\frac{1}{q_0}}\\
\leq C\biggl( \biggl(\int\limits_D \bigl|\nabla (u^{\frac{q}{q_0}})\bigr|^{p_0}w(x)~dx\biggr)^{\frac{1}{p_0}}+\biggl(\int\limits_D \bigl(|u|^{\frac{q}{q_0}}\bigr)^{p_0}w(x)~dx\biggr)^{\frac{1}{p_0}} \biggr)\\
\leq  C \biggl( \biggl(\int\limits_D |u|^{p_0\frac{q-q_0}{q_0}}|\nabla u|^{p_0}w(x)~dx\biggr)^{\frac{1}{p_0}}+\biggl(\int\limits_D |u|^{p_0\frac{q-q_0}{q_0}}|u|^{p_0}w(x)~dx\biggr)^{\frac{1}{p_0}} \biggr).
\nonumber
\end{multline}
Applying H\"older inequality we get
\begin{multline}
\biggl(\int\limits_D |u|^q w(x)~dx\biggr)^{\frac{1}{q_0}}
\leq C\biggl( \biggl(\int\limits_D |\nabla u|^{p}w(x)~dx\biggr)^{\frac{1}{p}}\cdot\biggl(\int\limits_D |u|^{pp_0\frac{q-q_0}{q_0(p-p_0)}}w(x)~dx\biggr)^{\frac{p-p_0}{pp_0}}\\
+\biggl(\int\limits_D |u|^{p}w(x)~dx\biggr)^{\frac{1}{p}}\cdot\biggl(\int\limits_D |u|^{pp_0\frac{q-q_0}{q_0(p-p_0)}}w(x)~dx\biggr)^{\frac{p-p_0}{pp_0}} \biggr).
\nonumber
\end{multline}
Because $q=pp_0\frac{q-q_0}{q_0(p-p_0)}$ we obtain finally
\begin{multline}
\biggl(\int\limits_D |u|^q w(x)~dx\biggr)^{\frac{1}{q_0}}\\
\leq C\biggl( \biggl(\int\limits_D |\nabla u|^{p}w(x)~dx\biggr)^{\frac{1}{p}}
+\biggl(\int\limits_D |u|^{p}w(x)~dx\biggr)^{\frac{1}{p}}\biggr)\cdot\biggl(\int\limits_D |u|^{q}w(x)~dx\biggr)^{\frac{p-p_0}{pp_0}}.
\nonumber
\end{multline}
Since $\frac{1}{q_0}-\frac{p-p_0}{pp_0}=\frac{1}{q}$ then
$$
\|u\mid L_q(D,w)\|\leq C\|u\mid W^1_p(D,w)\|.
$$
Lemma is proved.
\vskip 0.3cm

Denote 
$$
q^{\ast}_{m,D}(p)=\sup\{q\in \mathbb R^{+}:\,\,\text{the operator}\,\, W^m_p(D,w)\hookrightarrow L_q(D,w)\,\,\text{is bounded}\}.
$$

The following statements can be obtained directly from Lemma 3.
\vskip 0.3cm

{\bf Corollary~3.} If $D$ is a bounded domain in $\mathbb R^n$, $p\geq p_0$, then
$$
q^{\ast}_{1,D}(p) \geq \frac{pp_0 q^{\ast}_{1,D}(p_0)}{p_0 q^{\ast}_{1,D}(p_0)-p(q^{\ast}_{1,D}(p_0)-p_0)}.
$$
\vskip 0.3cm

{\bf Corollary~4.} If $D$ is a bounded domain in $\mathbb R^n$, $p\geq p_0$ and $m>1$, then
$$
q^{\ast}_{m,D}(p) \geq \frac{pp_0 q^{\ast}_{1,D}(p_0)}{p_0 q^{\ast}_{1,D}(p_0)-mp(q^{\ast}_{1,D}(p_0)-p_0)}.
$$
\vskip 0.3cm

{\sc Proof.} This corollary follows from the previous corollary by induction with respect to $m$.
\vskip 0.3cm

Combining Theorem~5 and  Corollary 4 we obtain finally. 

{\bf Theorem~9.} Let domain $D\subset\mathbb R^n$ be an embedding
domain and there exist $w$-weighted
$(p,q)$-quasiconformal homeomorphism $\varphi: D\to D'$ of domain $D$ onto bounded domain $D'$. 

If for some $s\leq
r\leq nq/(n-q)$ the following inequality is correct
$$
\int\limits_D\biggl(|J(x,\varphi)|w(\varphi(x))\biggr)^{\frac{r}{r-s}}<+\infty,
$$
then an embedding operator
$$
i: W^m_p(D',w)\hookrightarrow L_{s^{\ast}}(D',w)
$$
is bounded, if 
$$
s^{\ast}\geq \frac{ps}{s-m(s-p)}.
$$

\vskip 0.3cm

\bigskip

\centerline{\bf 5.~Solvability of degenerate elliptic equations}

\bigskip

In this section we apply an embedding theorem for Hilbert Sobolev spaces $W_2^1$ to a degenerate elliptic boundary problem.

{\bf Theorem~10.} Let $D'\subset\mathbb R^n$ be a quasiisometrical image of an embedding domain $D$ and
$$
\int\limits_{D'}w(y)^{-\frac{n}{2}}~dy<+\infty.
$$
Then an embedding operator 
$$
i' : W^1_2(D',w)\hookrightarrow L_2(D')
$$
is bounded.
\vskip 0.3cm

{\sc Proof.}
By conditions of the theorem there
exists a $Q$-quasiisometrical homeomorphism $\varphi : D\to D'$ of the embedding domain $D$ onto $D'$. For any
function $u\in W^1_2(D',w)$ the composition $u\circ \varphi$ is
weakly differentiable in domain $D$, and the following estimates are correct for any $1\leq q\leq 2$.
\begin{multline}
\|u\circ\varphi\mid W^1_q(D)\|=\biggl(\int\limits_D
|u\circ\varphi|^q~dx\biggr)^{\frac{1}{q}}+\biggl(\int\limits_D
|\nabla (u\circ\varphi)|^q~dx\biggr)^{\frac{1}{q}}\\
\leq
\biggl(\int\limits_D |u\circ\varphi|^q
(|J(x,\varphi)|w(\varphi(x))^{\frac{q}{2}}\frac{1}{(|J(x,\varphi)|w(\varphi(x))^{\frac{q}{2}}}
~dx\biggr)^{\frac{1}{q}}\\
+\biggl(\int\limits_D |\nabla u|^q
|\overline{\varphi}^{\prime}(x)|^q(|J(x,\varphi)|w(\varphi(x))^{\frac{q}{2}}\frac{1}{(|J(x,\varphi)|w(\varphi(x))^{\frac{q}{2}}}~dx\biggr)^{\frac{1}{q}}.
\nonumber
\end{multline}

By the H\"older inequality we have
\begin{multline}
\|u\circ\varphi\mid W^1_q(D)\|\leq\\
\leq
\biggl(\int\limits_D\biggl(\frac{1}{|J(x,\varphi)|w(\varphi(x))}\biggl)^{\frac{q}{2-q}}~dx\biggr)^{\frac{2-q}{2q}}
\biggl(\int\limits_D
|u|^2(\varphi(x))w(\varphi(x))|J(x,\varphi)|~dx\biggr)^{\frac{1}{2}}\\
+\biggl(\int\limits_D\biggl(\frac{|\overline{\varphi}^{\prime}(x)|^2}{|J(x,\varphi)|w(\varphi(x))}\biggl)^{\frac{q}{2-q}}~dx\biggr)^{\frac{2-q}{2q}}
\biggl(\int\limits_D |\nabla
u|^2(\varphi(x))w(\varphi(x))|J(x,\varphi)|~dx\biggr)^{\frac{1}{2}}.
\nonumber
\end{multline}

Since $\varphi$ is the $Q$-quasiisometrical homeomorphism, then
by the change of variable formula in the Lebesgue integral,
we obtain
\begin{multline}
\|u\circ\varphi\mid W^1_q(D)\|\leq Q^{\frac{n}{q}}
\biggl(\int\limits_{D'}w(y)^{\frac{q}{q-2}}~dy\biggr)^{\frac{2-q}{2q}}\biggl(\int\limits_{D'}|u|^2
w(y) ~dy\biggr)^{\frac{1}{2}}\\
+Q^{\frac{2-q+n}{q}}
\biggl(\int\limits_{D'}w(y)^{\frac{q}{q-2}}~dy\biggr)^{\frac{2-q}{2q}}\biggl(\int\limits_{D'}|\nabla
u|^p w(y)~dy\biggr)^{\frac{1}{2}}\\
=Q^{\frac{n}{q}}K(w)\|u\mid
L_2(D',w)\|+Q^{\frac{2-q+n}{q}}K(w)\|\nabla u\mid L_2(D'w)\|.
\nonumber
\end{multline}

By the previous inequality the composition operator 
$$
\varphi^{\ast}: W^1_2(D',w)\to W^1_q(D),\,\,\, 1\leq q\leq
2,
$$
is bounded.

Since $D$ is an embedding domain there exists the bounded embedding operator
$$
i: W^1_q(D)\hookrightarrow L_{\frac{nq}{n-q}}(D).
$$
Now we choose a number $q$ such that $\frac{nq}{n-q}=2$ (i.e. $q=\frac{2n}{n+2}$).

Since $\varphi$ is the $Q$-quasiisometrical homeomorphism,
the following composition operator acting on Lebesgue spaces
$$
(\varphi^{-1})^{\ast}: L_2(D)\to L_2(D'),
$$
is bounded also [18].

Therefore the embedding operator $i': W^1_2(D',w)\hookrightarrow
L_2(D')$ is bounded as a composition of bounded operators
$\varphi^{\ast}$, $(\varphi^{-1})^{\ast}$ and the bounded embedding
operator $i$.

\vskip 0.3cm

Define an inner product in the weighted space $W^1_2(D,w)$ as:
$$
<u,v>=\int\limits_D \nabla u\cdot\nabla v ~w(x)~dx,
$$
for any $u,v\in W^1_2(D,w)$.

\vskip 0.3cm

Consider Dirichlet problem for the degenerate elliptic equation:
\begin{gather}
\dv(w(x)\nabla u)=f,\\
u\vert_{\partial D}=0.
\end{gather}
in a bounded domain $D$ for the weight $w\in C^1(D)$.

\vskip 0.3cm

{\bf Theorem~11.} Let $f\in L_2(D)$ and $\int\limits_{D}w(y)^{-\frac{n}{2}}~dx<+\infty$. Then there exists the unique weak solution $u\in \overset{\circ}{W}^1_2(D,w)$ of the problem $(1,2)$.
\vskip 0.3cm

{\sc Proof.}
The function $f\in L_2(D)$ induces a linear functional $F: L_2(D)\rightarrow R$ by the standard rule
$$
F(\phi)=\int\limits_D f(x)\phi(x)~dx.
$$
By Theorem~10 there exists bounded embedding operator $i: \overset{\circ}{W}^1_2(D,w)\hookrightarrow L_2(D)$. 

Therefore 
\begin{multline}
|F(\phi)|\leq\|f\cdot\phi\mid L_1(D)\|\\
\leq\|f\mid L_2(D)\|\cdot\|\phi\mid L_2(D)\|\leq C\|f\mid L_2(D)\|\cdot\|\phi\mid \overset{\circ}{W}^1_2(D,w)\|.
\nonumber
\end{multline}
Hence, $F$ is a bounded linear functional in Hilbert space $\overset{\circ}{W}^1_2(D,w)$. By Riesz representation theorem [20] there exists the unique function $u\in \overset{\circ}{W}^1_2(D,w)$ such that

$$
F(\phi)=<u,\phi>=\int\limits_D \nabla u\cdot\nabla\phi~w(x)dx,
$$
or
$$
\int\limits_D \nabla u\cdot\nabla\phi~w(x)dx=\int\limits_D f(x)\phi(x)~dx.
$$
Therefore $u$ is the unique weak solution of the problem $(1,2)$.

\vskip 0.5cm

\centerline{REFERENCES}

\begin{enumerate}
\item
Kufner A. {\it Weighted Sobolev spaces.} -- Leipzig: 
Teubner-Texte fur Mathematik, 1980.

\item
Turesson B.~O. {\it Nonlinear Potential Theory and Weighted Sobolev spaces.} -- Berlin: 
Springer. (Lecture Notes in Mathematics), 2000.

\item
Gurka P., Opic B. {\it Continuous and compact embeddings of weighted Sobolev spaces. I.}//
Czech. Math. J. -- 1988. -- V.~38. -- N.~113. -- P.~730--744.
 
\item
Gurka P., Opic B. {\it Continuous and compact embeddings of weighted Sobolev spaces. II.}//
Czech. Math. J. -- 1989. -- V.~39. -- N.~114. -- P.~78--94.
 
\item
Gurka P., Opic B. {\it Continuous and compact embeddings of weighted Sobolev spaces. III.}//
Czech. Math. J. -- 1991. -- V.~41. -- N.~116. -- P.~317--341.

\item
Antoci F. {\it Some necessary and some sufficient conditions for the compactness of the embedding of weighted Sobolev spaces.}// Ricerche Mat. -- 2003. -- V.~52. -- N.~1. -- P.~55--71.

\item
Besov O.~V. {\it The Sobolev embedding theorem for a domain with an irregular boundary.}//
Doklady Mathematics -- 2000. -- V.~62. -- N.~ 1. -- P.~22--25.

\item
David G., Semmes S. {\it Strong $A_{\infty}$-weights, Sobolev inequalities and quasiconformal mappings.}//
Lecture notes in Pure and Appl. Math. -- 1990. -- V.~122. --  P.~101--111.

\item
Gol'dshtein V., Gurov L. {\it Applications of change of variables operators for exact embedding theorems.}//
Integral Equations Operator Theory -- 1994. -- V.~19. -- N.~1. -- P.~1--24.

\item
Gol'dshtein V., Ramm A.~G. {\it Compactness of the embedding operators for rough domains.}//
Math. Inequalities and Applications -- 2001. -- V.~4. -- N.~1. -- P.~127--141.
 
\item
Vodop'yanov S.~K., Ukhlov A. {\it Set functions and their applications in the theory of Lebesgue and Sobolev spaces. II}//
Siberian Adv. Math. -- 2005. -- V.~15. -- N.~1. -- P.~91--125.

\item
Burenkov V.~I. {\it Sobolev Spaces on Domains}- Stuttgart: Teubner-Texter zur Mathematik. 1998.

\item
Stein E. {\it Harmonic analysis real-variable methods, orthogonality, and oscillatory integrals}--
Princeton: Princeton Univ. Press. 1993.

\item
Kilpelainen~T. {\it Weighted Sobolev spaces and capacity.}//
Ann. Acad. Sci. Fenn. Ser. A. I. Math. -- 1994. -- V.~19. -- P.~95--113.

\item
Heinonen J., Kilpelainen T., Martio O. {\it Nonlinear potential theory of degenerate elliptic equations}--
Oxford: Oxford Univ. Press. 1993.

\item
Maz'ya V. {\it Sobolev spaces} -- Berlin: Springer Verlag. 1985.

\item
Gol'dshtein V.~M., Reshetnyak Yu.~G. {\it Quasiconformal mappings and Sobolev spaces} -- Dordrecht, Boston, London: Kluwer Academic Publishers. 1990.

\item
Vodop'yanov S.~K., Ukhlov A. {\it Set functions and their applications in the theory of Lebesgue and Sobolev spaces. I}//
Siberian Adv. Math. -- 2004. -- V.~14. -- N.~4. -- P.~78--125.

\item
Vodop'yanov S.~K., Ukhlov A. {\it Mappings associated with weighted Sobolev spaces.}//
Contemporary Mathematics  (to appear)

\item
Riesz F. {\it Sur les operations fonctionelles lineaires.}//
C. R. Acad. Sci. Paris  -- 1909. -- V.~149. -- P.~974--977.

\end{enumerate}

\end{document}